\documentclass[11pt,leqno]{article}
\usepackage[margin=1in]{geometry} 
\usepackage{amssymb,amsfonts,amsmath,bbm,mathrsfs,stmaryrd,mathtools}
\usepackage{xcolor}
\usepackage{url}

\usepackage{extarrows}

\usepackage[shortlabels]{enumitem}
\usepackage{tensor}

\usepackage{xr}
\usepackage[T1]{fontenc}

\usepackage[colorlinks,
linkcolor=black!75!red,
citecolor=blue,
pdftitle={},
pdfproducer={pdfLaTeX},
pdfpagemode=None,
bookmarksopen=true,
bookmarksnumbered=true,
backref=page]{hyperref}

\usepackage{tikz}
\usetikzlibrary{arrows,calc,decorations.pathreplacing,decorations.markings,intersections,shapes.geometric,through,fit,shapes.symbols,positioning,decorations.pathmorphing}

\makeatletter
\newcommand*\circled[1]{\tikz[baseline=(char.base)]{
    \node[shape=circle, draw, inner sep=0pt, 
    minimum height={\f@size},] (char) {\vphantom{WAH1g}#1};}}
\makeatother
    
\usepackage{braket}

\usepackage[amsmath,thmmarks,hyperref]{ntheorem}
\usepackage{cleveref}

\newcommand\nthalias[1]{\AddToHook{env/#1/begin}{\crefalias{lemma}{#1}}}

\nthalias{definition}
\nthalias{example}
\nthalias{examples}
\nthalias{remark}
\nthalias{remarks}
\nthalias{convention}
\nthalias{notation}
\nthalias{construction}
\nthalias{theoremN}
\nthalias{propositionN}
\nthalias{corollaryN}
\nthalias{lemma}
\nthalias{proposition}
\nthalias{corollary}
\nthalias{theorem}
\nthalias{conjecture}
\nthalias{question}
\nthalias{assumption}

\creflabelformat{enumi}{#2#1#3}

\crefname{section}{Section}{Sections}
\crefformat{section}{#2Section~#1#3} 
\Crefformat{section}{#2Section~#1#3} 

\crefname{subsection}{\S}{\S\S}
\AtBeginDocument{%
  \crefformat{subsection}{#2\S#1#3}%
  \Crefformat{subsection}{#2\S#1#3}%
}

\crefname{subsubsection}{\S}{\S\S}
\AtBeginDocument{%
  \crefformat{subsubsection}{#2\S#1#3}%
  \Crefformat{subsubsection}{#2\S#1#3}%
}

\theoremstyle{plain}

\newtheorem{lemma}{Lemma}[section]
\newtheorem{proposition}[lemma]{Proposition}
\newtheorem{corollary}[lemma]{Corollary}
\newtheorem{theorem}[lemma]{Theorem}

\theoremstyle{plain}
\theoremnumbering{Alph}
\newtheorem{theoremN}{Theorem}

\theoremstyle{plain}
\theorembodyfont{\upshape}
\theoremsymbol{\ensuremath{\blacklozenge}}

\newtheorem{definition}[lemma]{Definition}
\newtheorem{example}[lemma]{Example}

\newtheorem{remark}[lemma]{Remark}
\newtheorem{remarks}[lemma]{Remarks}

\newtheorem{notation}[lemma]{Notation}

\crefname{definition}{definition}{definitions}
\crefformat{definition}{#2definition~#1#3} 
\Crefformat{definition}{#2Definition~#1#3} 

\crefname{ex}{example}{examples}
\crefformat{example}{#2example~#1#3} 
\Crefformat{example}{#2Example~#1#3} 

\crefname{exs}{example}{examples}
\crefformat{examples}{#2example~#1#3} 
\Crefformat{examples}{#2Example~#1#3} 

\crefname{remark}{remark}{remarks}
\crefformat{remark}{#2remark~#1#3} 
\Crefformat{remark}{#2Remark~#1#3} 

\crefname{remarks}{remark}{remarks}
\crefformat{remarks}{#2remark~#1#3} 
\Crefformat{remarks}{#2Remark~#1#3} 

\crefname{convention}{convention}{conventions}
\crefformat{convention}{#2convention~#1#3} 
\Crefformat{convention}{#2Convention~#1#3} 

\crefname{notation}{notation}{notations}
\crefformat{notation}{#2notation~#1#3} 
\Crefformat{notation}{#2Notation~#1#3} 

\crefname{table}{table}{tables}
\crefformat{table}{#2table~#1#3} 
\Crefformat{table}{#2Table~#1#3}

\crefname{lemma}{lemma}{lemmas}
\crefformat{lemma}{#2lemma~#1#3} 
\Crefformat{lemma}{#2Lemma~#1#3} 

\crefname{proposition}{proposition}{propositions}
\crefformat{proposition}{#2proposition~#1#3} 
\Crefformat{proposition}{#2Proposition~#1#3} 

\crefname{propositionN}{proposition}{propositions}
\crefformat{propositionN}{#2proposition~#1#3} 
\Crefformat{propositionN}{#2Proposition~#1#3} 

\crefname{corollary}{corollary}{corollaries}
\crefformat{corollary}{#2corollary~#1#3} 
\Crefformat{corollary}{#2Corollary~#1#3} 

\crefname{corollaryN}{corollary}{corollaries}
\crefformat{corollaryN}{#2corollary~#1#3} 
\Crefformat{corollaryN}{#2Corollary~#1#3} 

\crefname{theorem}{theorem}{theorems}
\crefformat{theorem}{#2theorem~#1#3} 
\Crefformat{theorem}{#2Theorem~#1#3} 

\crefname{theoremN}{theorem}{theorems}
\crefformat{theoremN}{#2theorem~#1#3} 
\Crefformat{theoremN}{#2Theorem~#1#3} 

\crefname{enumi}{}{}
\crefformat{enumi}{#2#1#3}
\Crefformat{enumi}{#2#1#3}

\crefname{assumption}{assumption}{Assumptions}
\crefformat{assumption}{#2assumption~#1#3} 
\Crefformat{assumption}{#2Assumption~#1#3} 

\crefname{construction}{construction}{Constructions}
\crefformat{construction}{#2construction~#1#3} 
\Crefformat{construction}{#2Construction~#1#3} 

\crefname{question}{question}{Questions}
\crefformat{question}{#2question~#1#3} 
\Crefformat{question}{#2Question~#1#3} 

\crefname{equation}{}{}
\crefformat{equation}{(#2#1#3)} 
\Crefformat{equation}{(#2#1#3)}

\numberwithin{equation}{section}

\theoremstyle{nonumberplain}
\theoremsymbol{\ensuremath{\blacksquare}}

\newtheorem{proof}{Proof}
\newcommand\pf[1]{\newtheorem{#1}{Proof of \Cref{#1}}}

\newcommand\bC{{\mathbb C}}

\newcommand\bG{{\mathbb G}}

\newcommand\bN{{\mathbb N}}

\newcommand\bP{{\mathbb P}}

\newcommand\bR{{\mathbb R}}
\newcommand\bS{{\mathbb S}}

\newcommand\bZ{{\mathbb Z}}

\newcommand\cC{{\mathcal C}}

\newcommand\cP{{\mathcal P}}

\newcommand\cT{{\mathcal T}}

\newcommand\wt{\widetilde}

\DeclareMathOperator{\Ad}{Ad}
\DeclareMathOperator{\id}{id}

\DeclareMathOperator{\im}{\mathrm{im}}

\DeclareMathOperator{\spn}{\mathrm{span}}

\DeclareMathOperator{\GL}{GL}
\DeclareMathOperator{\PGL}{PGL}

\DeclareMathOperator{\SU}{SU}

\DeclareMathOperator{\D}{D}
\DeclareMathOperator{\G}{G}
\DeclareMathOperator{\T}{T}
\DeclareMathOperator{\U}{U}

\DeclareMathOperator{\diag}{diag}

\newcommand{\qedhere}{\mbox{}\hfill\ensuremath{\blacksquare}}

\newcommand{\xrightarrowdbl}[2][]{%
  \xrightarrow[#1]{#2}\mathrel{\mkern-14mu}\rightarrow
}

\title{Spectral selections, commutativity preservation and Coxeter-Lipschitz maps}
\author{Alexandru Chirvasitu}

\begin{document}

\date{}

\newcommand{\Addresses}{{
  \bigskip
  \footnotesize

  \textsc{Department of Mathematics, University at Buffalo}
  \par\nopagebreak
  \textsc{Buffalo, NY 14260-2900, USA}  
  \par\nopagebreak
  \textit{E-mail address}: \texttt{achirvas@buffalo.edu}

}}

\maketitle

\begin{abstract}
  For a Coxeter system $(W,S)$ a self-map $\tau$ of $W$ such that $\tau_{\sigma\theta}\in \{\tau_{\theta},\ \sigma\tau_{\theta}\}$ for all $\theta\in W$ and all reflections $\sigma$ (analogous to being 1-Lipschitz with respect to the Bruhat order on $W$) is either constant or a right translation. A somewhat stronger version holds for $S_n$, where it suffices that $\sigma$ range over smaller, $\theta$-dependent sets of reflections.
  
  These combinatorial results have a number of consequences concerning continuous spectrum- and commutativity-preserving maps $\mathrm{SU}(n)\to M_n$ defined on special unitary groups: every such map is a conjugation composed with (a) the identity; (b) transposition, or (c) a continuous diagonal spectrum selection. This parallels and recovers Petek's analogous statement for self-maps of the space $H_n\le M_n$ of self-adjoint matrices, strengthening it slightly by expanding the codomain to $M_n$. 
\end{abstract}


\vspace{.5cm}

\noindent{MSC 2020: 20B30; 20F55; 54C05; 47A10; 20F10; 15B57; 06A06; 54H15
  
  
}


\section*{Introduction}

Among the few strands of motivation for the material below are the various classification results for maps between matrix spaces preserving spectral and/or algebraic structure, variations of which abound in the literature: \cite{MR2853129,2501.06840v2,2501.18776v1,GogicPetekTomasevic,GOGIC2025129497,MR832991,Petek-HM,Petek-TM,zbMATH01100760,PetekSemrl,MR859735,MR1311919,Semrl2,Semrl} and \emph{their} references will provide a still-small sample. Consider, as one concrete entry point, the main result of \cite{Petek-HM}. Here and throughout the paper we refer to a(n often partial) self-function $M_n\xrightarrow{\phi}M_n$ of the space of $n\times n$ matrices as
\begin{itemize}
\item \emph{commutativity-}(or just \emph{C-})\emph{preserving} if $\phi(X)$ and $\phi(Y)$ commute whenever $X$ and $Y$ to;
\item \emph{spectrum-}(or just \emph{S-})\emph{preserving} if $X$ and $\phi(X)$ have the same spectrum (as a subset of $\bC^n$);

\item and \emph{CS-preserving} if both conditions are met.
\end{itemize}
Denote by $\Ad_g$ the conjugation action $g\cdot g^{-1}$ on a group $\G$ by an element $g\in \G$. \cite[Main theorem]{Petek-HM}, then, says that for positive integers $n\ge 3$ the continuous, CS-preserving self-maps of the (real) algebra $H_n\le M_n:=M_n(\bC)$ of $n\times n$ self-adjoint matrices are precisely those of the form
\begin{itemize}[wide]
\item $\Ad_T (-)^{\bullet}$ for some unitary $T\in \U(n)$ where the symbol $\bullet$ is either blank or `$t$', denoting transposition;
\item or
  \begin{equation}\label{eq:sa.reord.spc}
    X
    \xmapsto{\quad}
    \Ad_T \diag\left(\eta_j(X),\ 1\le j\le n\right)
    ,\quad
    T\in \U(n)
  \end{equation}
  where
  \begin{equation*}
    \eta_1(X)\le \cdots\le \eta_n(X)
  \end{equation*}
  is the non-increasing ordering of the (real) spectrum of $X=X^*$. 
\end{itemize}

The present note is partly concerned with a variant of that result valid for special unitary matrices instead. 

\begin{theoremN}\label{thn:cls.sun}
  Let $n\in \bZ_{\ge 3}$. The continuous CS-preserving maps $\SU(n)\to M_n$ are precisely those of the form
  \begin{enumerate}[(a),wide]
  \item\label{item:thn:cls.sun:cj} $\Ad_T (-)^{\bullet}$ for $\bullet\in \{\text{blank},\ t\}$ and some $T\in \GL(n)$;

  \item\label{item:thn:cls.sun:ord} or
    \begin{equation*}
      \D
      \ni
      X
      \xmapsto{\quad}
      \Ad_T(\lambda_j(X))_j
      \in M_n
      \quad
      \left(\text{some $T\in \GL(n)$}\right),
    \end{equation*}
    where $\lambda_j(X):=\exp(2\pi i x_j)$ for the unique
    \begin{equation*}
      x_1\le x_2\le \cdots\le x_n\le x_1+1
      ,\quad
      \sum_j x_j=0
    \end{equation*}
    for which $\exp(2\pi i x_j)$ constitute the spectrum of $X$.  \qedhere
  \end{enumerate}
\end{theoremN}

Item \Cref{item:thn:cls.sun:ord} in \Cref{thn:cls.sun} (where tuples are meant as the respective diagonal matrices) perhaps requires some unpacking. It follows from \cite[Lemmas 1 and 2]{MR210096}, auxiliary to describing the $n^{th}$ symmetric power
\begin{equation*}
  \left(\bS^1\right)^{[n]}
  :=
  \left(\bS^1\right)^n/S_n
  ,\quad
  S_n:=\text{$n$-symbol permutation group}
\end{equation*}
of the circle, that
\begin{equation}\label{eq:morton.ident}
  \left\{(x_j)\in \bR^n
    \ :\
    \begin{aligned}
      x_1\le x_2\le &\cdots\le x_n\le x_1+1\\
      \sum_j x_j&=0
    \end{aligned}
  \right\}
  \xmapsto[\quad\cong\quad]{\quad}
  \left(\exp 2\pi i x_j\right)_j
\end{equation}
is (after identifying tuples up to permutation) a homeomorphism onto
\begin{equation}\label{eq:tor.sun}
  \left\{(\zeta_j)\in \left(\bS^1\right)^{[n]}\ :\ \prod_j \zeta_j=1\right\}
  \subset
  \left(\bS^1\right)^{[n]}. 
\end{equation}
This, incidentally, is intimately linked to the geometry of \emph{Weyl-group} actions on \emph{maximal tori} \cite[Chapter IV]{btd_lie_1995} of compact Lie groups: \Cref{eq:tor.sun} is the quotient $\T/W$ of that action on the maximal torus
\begin{equation*}
  \T:=\left\{(\zeta_j)\in \left(\bS^1\right)^{n}\ :\ \prod_j \zeta_j=1\right\}
  \le
  \SU(n)
  ,\quad
  W=S_n\text{ acting by permutations},
\end{equation*}
and Morton's homeomorphism \Cref{eq:morton.ident} is an instance of the fact \cite[\S 4.8, Theorem]{hmph_cox} that for a simply-connected, simple compact Lie group $\G$ the respective quotient
\begin{equation*}
  \G/\left(\text{adjoint action}\right)
  \xrightarrow[\quad\cong\quad]{\quad\text{\cite[Proposition IV.2.6]{btd_lie_1995}}\quad}
  \T/W
\end{equation*}
is always identifiable with a simplex in the Lie algebras $Lie(T)$. This is what underlies the continuous \emph{eigenvalue selection}
\begin{equation*}
  \SU(n)\ni X
  \xmapsto{\quad}
  \lambda_1(X)
  \quad
  \left(\text{notation of \Cref{thn:cls.sun}}\right)
\end{equation*}
in \cite[Remark 1.4(3)]{2501.06840v2}, applicable more generally \cite[Theorem A]{2502.08847v1} to simply-connected compact Lie groups. As the preceding discussion makes clear,
\begin{equation*}
    \SU(n)\ni X
    \xmapsto{\quad}
    \left(\lambda_j(X)\right)_{j}
    \in
    \left(\bS^1\right)^n
\end{equation*}
is a continuous spectrum ordering (the title's \emph{spectral selection}). 

Essentially the same combinatorial principle (\Cref{pr:chartaus}) ultimately driving \Cref{thn:cls.sun} will also recover \cite[Main theorem]{Petek-HM} (in a slightly stronger form: the codomain is all of $M_n$ as opposed to $H_n$). 

\begin{theoremN}\label{thn:cls.sa}
  Let $n\in \bZ_{\ge 3}$. The continuous CS-preserving maps $H_n\to M_n$ are precisely those of the form
  \begin{enumerate}[(a),wide]
  \item\label{item:thn:cls.sa:cj} $\Ad_T (-)^{\bullet}$ for $\bullet\in \{\text{blank},\ t\}$ and some $T\in \GL(n)$;

  \item\label{item:thn:cls.sa:ord} or of the form \Cref{eq:sa.reord.spc} for some $T\in \GL(n)$.  \qedhere
  \end{enumerate}
\end{theoremN}

The combinatorial content of \Cref{thn:cls.sun,thn:cls.sa} in turn suggests and motivates the offshoot material of \Cref{se:cox.lip}. An examination of the proof of \Cref{thn:cls.sun} via \Cref{pr:chartaus} below (in parallel to that of \cite[Theorem 2.1]{2501.06840v2} by means of \cite[Proposition 2.2]{2501.06840v2}) makes it clear that the constraints imposed by spectrum and commutativity preservation on maps defined on
\begin{itemize}[wide]
\item maximal tori (in the case of $\SU(n)$ or $\U(n)$);

\item or maximal abelian subalgebras (in the case of $H_n$ pertinent to \Cref{thn:cls.sa})
\end{itemize}
are intimately connected to the metric geometry of the \emph{Coxeter complex} (\cite[\S 2.1]{ron_build_2e_2009}, \cite[Chapter 3, Exercise 16]{bb_comb-cox}) attached to the usual \cite[Example 1.2.3]{bb_comb-cox} realization
\begin{equation*}
  \left(S_n,\ \left\{\text{transpositions }(j\ j+1)\right\}_{j=1}^{n-1}\right)
\end{equation*}
of the symmetric group on $n$ symbols as (the underlying group of) a \emph{Coxeter system} \cite[\S 1.1, p.2]{bb_comb-cox}.  

A trimmed-down, paraphrased aggregate of \Cref{pr:chartaus} and \Cref{th:cox.lip} below, stemming ultimately from an examination of the combinatorics of diagonally-defined CS-preserving maps, reads as follows (with $\Ad_W S$, for a subset $S\subseteq W$ of a group, denoting the set of $W$-conjugates of $S$-elements).

\begin{theoremN}\label{thn:cx}
  Let $(W,S)$ be a Coxeter system.

  \begin{enumerate}[(1),wide]
  \item A self-map $W\xrightarrow{\tau}W$ satisfying
    \begin{equation}\label{eq:lip.cond.intro}
      \forall\left(\theta\in W\right)
      \forall\left(\sigma\in \Ad_W S\in W\right)
      \quad:\quad
      \tau_{\sigma\theta}\in \left\{\tau_{\theta},\ \sigma\tau_{\theta}\right\}
    \end{equation}
    is either constant or a right translation. 

  \item The same conclusion holds in the symmetric-group case
    \begin{equation*}
      (W,S)
      :=
      (S_n,\ \left\{(1\ 2),\ \cdots,\ (n-1\ n)\right\})
    \end{equation*}
    if \Cref{eq:lip.cond.intro} is assumed only for $\theta\in \Ad_{\theta}\big(S\sqcup\{(n\ 1)\}\big)$.  \qedhere
  \end{enumerate}
\end{theoremN}

\subsection*{Acknowledgments}

I am indebted to I. Gogi{\'c}, J. Mangahas and M. Toma{\v{s}}evi{\'c} for numerous insightful comments, as well as the anonymous referees' suggestions and pointers. 


\section{Circle configuration spaces and special unitary groups}\label{se:sun}

The \emph{configuration space} \cite[p.vii]{fh_config} $\cC^n(\bS^1)$ of $n$-tuples in $\bS^1$ with distinct entries can be thought of as the space of simple-spectrum diagonal unitary matrices. It is not unnatural to pose the analogous problem for the subspace \emph{special} unitary matrices, and to note the distinction between the two cases.

Recall that a group acts \emph{simply transitively} on a set (\emph{sharply 1-transitively} in \cite[Definition post Theorem 9.7]{rot-gp}) if the action is both free and transitive. 

\begin{lemma}\label{le:su.diag.1trans}
  For $n\in \bN$ the action of the symmetric group $S_n$ on the connected (path-)components of
  \begin{equation*}
    \cC^n(\bS^1)_{\Pi=1}
    :=
    \left\{(z_i)\in \cC^n(\bS^1)\ :\ \prod_i z_i=1\right\}
  \end{equation*}
  is simply transitive. 
\end{lemma}
\begin{proof}    
  The discrete group $S_n$ acts freely on the manifold $\cC^n(\bS^1)_{\Pi=1}$, so
  \begin{equation}\label{eq:princ.sn.fib}
    \cC^n(\bS^1)_{\Pi=1}
    \xrightarrowdbl{\quad}
    \cC^n(\bS^1)_{\Pi=1}/S_n
  \end{equation}
  is an \emph{$S_n$-principal covering} in the sense of \cite[\S 14.1]{td_alg-top}. By \cite[Lemmas 1 and 2]{MR210096} the base space $\cC^n(\bS^1)_{\Pi=1}/S_n$ of that covering can be identified with the interior
  \begin{equation*}
    \left\{(x_i)_{i=1}^n\in \bR^n\ :\ \sum_i x_i=0\quad\text{and}\quad x_1<x_2<\cdots<x_n<x_1+1\right\}
  \end{equation*}
  of a simplex, and is thus \emph{contractible} \cite[\S 2.1]{td_alg-top}. It follows \cite[Theorem 14.4.1]{td_alg-top} that the principal fibration \Cref{eq:princ.sn.fib} is trivial:
  \begin{equation*}
    \cC^n(\bS^1)_{\Pi=1}
    \cong
    (\cC^n(\bS^1)_{\Pi=1}/S_n)\times S_n
  \end{equation*}
  as $S_n$-spaces, concluding the proof. 
\end{proof}

The following result is a special unitary version of \cite[Proposition 2.2]{2501.06840v2}. 

\begin{theorem}\label{th:su.tor}
  For $n\in \bZ_{>0}$ a continuous CS-preserving map $\D\xrightarrow{\phi}M_n$ defined on the $n\times n$ diagonal special unitary group $\D\le \SU(n)$ is either
  \begin{enumerate}[(a),wide]
  \item\label{item:th:su.tor:cj} conjugation $\Ad_T$ by some $T\in \GL(n)$;

  \item\label{item:th:su.tor:ord} or
    \begin{equation*}
      \D
      \ni
      X
      \xmapsto{\quad}
      \Ad_T(\lambda_j(X))_j
      \in M_n,
    \end{equation*}
    where $\lambda_j(X):=\exp(2\pi i x_j)$ for the unique
    \begin{equation*}
      x_1\le x_2\le \cdots\le x_n\le x_1+1
      ,\quad
      \sum_j x_j=0
    \end{equation*}
    for which $\exp(2\pi i x_j)$ constitute the spectrum of $X$.
  \end{enumerate}  
\end{theorem}
\begin{proof}
  There is a common core to the present argument and that proving the unitary branch of \cite[Proposition 2.2]{2501.06840v2}. First, there is no loss in assuming $\phi$ takes diagonal values: this is so after composing with a conjugation for simple-spectrum matrices, and hence also generally by continuity. This means that $\phi$ simply permutes diagonal entries:
  \begin{equation*}
    \forall\left(Y=(y_j)_j\in \D\right)
    \exists\left(\tau\in S_n\right)
    \quad:\quad
    \phi(Y)=(y_{\tau j})_j.
  \end{equation*}
  The permutation $\tau$ is constant along (path-)components of $\D$, and uniquely determined for simple-spectrum $Y$. This gives a map
  \begin{equation}\label{eq:comps2sn}
    \pi_0\left(\cC^n(\bS^1)_{\Pi=1}\right)
    \ni
    C
    \xmapsto{\quad}
    \tau_C\in S_n,
  \end{equation}
  and hence also a self-map $\tau_{\bullet}$ on $S_n$ after identifying the domain $\pi_0\left(\cC^n(\bS^1)_{\Pi=1}\right)$ of \Cref{eq:comps2sn} with $S_n$ equivariantly, as allowed by \Cref{le:su.diag.1trans}. For the same reasons as in the proof of \cite[Proposition 2.2]{2501.06840v2}, we have
  \begin{equation}\label{eq:trnsp.lipsch}
    \forall\left(\theta\in S_n\right)
    \forall\left(\text{c-simple transposition }\sigma\in S_n\right)      
    \quad:\quad
    \tau_{\theta\sigma}
    =
    \tau_{\Ad_{\theta}\sigma\cdot \theta}
    \in
    \left\{
      \tau_{\theta}
      ,\
      \Ad_{\theta}\sigma\cdot \tau_{\theta}
    \right\},
  \end{equation}
  where
  \begin{equation*}
    \left\{
      \text{\emph{c-simple} transpositions (`c' for `cyclic')}
    \right\}
    :=
    \left\{
      (1\ 2)
      ,\ 
      (2\ 3)
      ,\ 
      \cdots
      ,\ 
      (n-1\ n)
      ,\ 
      (n\ 1)
    \right\}.
  \end{equation*}  
  By \Cref{pr:chartaus} $\tau_{\bullet}$ is either constant or right translation by some element of $S_n$; the two options respectively corresponding to those of the statement, the proof is complete.
\end{proof}

Note also the parallels to the description of continuous CS-preserving self-maps of the space $H_n$ of Hermitian $n\times n$ matrices in \cite[Main theorem]{Petek-HM}. The immediately-guessable Hermitian analogue of \Cref{th:su.tor}, however, does \emph{not} hold: \Cref{ex:herm.hyb} below provides a continuous self-map of the diagonal Hermitian matrices which is neither a conjugation $\Ad_T$, $T\in \GL(n)$ nor of the form
\begin{equation*}
  X
  \xmapsto{\quad}
  \Ad_T(\lambda_j(X))_j
  \in M_n,
\end{equation*}
for the non-decreasing enumeration
\begin{equation*}
  \lambda_1(X)\le \lambda_2(X)\le \cdots\le \lambda_n(X)
\end{equation*}
of the spectrum of $X$.  

\begin{example}\label{ex:herm.hyb}
  We describe a self-map of the space $H_{n,d}$ of diagonal $n\times n$ Hermitian matrices (i.e. diagonal and real) for $n=3$ (the `$d$' subscript is for `diagonal'). We employ the notation
  \begin{equation*}
    [c\ a\ b]
    :=
    \begin{pmatrix}
      c&0&0\\
      0&a&0\\
      0&0&b
    \end{pmatrix}
    ,\quad \left(a\le b\le c\right)\subset \bR.
  \end{equation*}
  With those conventions, the map in question will be
  \begin{alignat*}{2}
    [a\ b\ c]&\xmapsto{\quad}[a\ b\ c]
               ,\quad
               [a\ c\ b]&\xmapsto{\quad} [a\ c\ b]\\
    [b\ a\ c]&\xmapsto{\quad}[a\ b\ c]
               ,\quad
               [b\ c\ a]&\xmapsto{\quad} [a\ c\ b]\\
    [c\ a\ b]&\xmapsto{\quad}[a\ b\ c]
               ,\quad
               [c\ b\ a]&\xmapsto{\quad} [a\ b\ c]\\
  \end{alignat*}
  The (easily checked) claim is that this is indeed a well-defined continuous self-map of $H_{3,d}$: there are compatibility constraints for the cases $a=b$ or $b=c$ (or both), and all such hold. 
\end{example}


\begin{proposition}\label{pr:chartaus}
  For $n\in \bZ_{\ge 1}$ the self-maps $\tau_{\bullet}$ of the symmetric group $S_n$ satisfying \Cref{eq:trnsp.lipsch} are precisely the constants and the right translations. 
\end{proposition}
\begin{proof}
  Some language, in light of \Cref{eq:trnsp.lipsch}:
  \begin{itemize}[wide]
  \item $\theta\in S_n$ \emph{($\tau$-)grows along $\sigma$} if $\tau_{\theta\sigma} = \Ad_{\theta}\sigma\cdot \tau_{\theta}$;

  \item $\theta\in S_n$ \emph{($\tau$-)lingers along $\sigma$} if $\tau_{\theta\sigma} = \tau_{\theta}$
  \end{itemize}
  for c-simple $\sigma$. I claim that either
  \begin{enumerate}[(a),wide]
  \item\label{item:pr:chartaus:grw} all $\theta\in S_n$ grow along all c-simple $\sigma$, or
  \item\label{item:pr:chartaus:lngr} all $\theta\in S_n$ linger along all c-simple $\sigma$.
  \end{enumerate}
  That this then completes the proof is clear: \Cref{item:pr:chartaus:grw} renders $\tau$ a right translation, while \Cref{item:pr:chartaus:lngr} means it must be constant. The rest of the proof is thus devoted to the claim itself.
  
  Some notation will aid the argument. First, denote the c-simple transpositions by
  \begin{equation*}
    \xi_j:=(j\ j+1)
    ,\quad
    j\in [n]:=\left\{1\cdots n\right\}
    \quad\text{with}\quad
    n+1=1
    \quad\left(\text{and hence }\xi_n=(n\ 1)\right).
  \end{equation*}
  Next, for any $\theta\in S_n$ and $j\in [n]$ we define $\theta_{j\bullet}$ recursively by
  \begin{equation*}
    \forall \left(k\in [n]\right)
    \quad:\quad
    \theta_{jk}\cdots \theta_{j2}\cdot \theta_{j1}\cdot \theta
    =
    \theta\cdot\xi_j\cdot \xi_{j+1}\cdots \xi_{j+k-1}. 
  \end{equation*}
  Explicitly:
  \begin{equation*}
    \theta_{jk}
    =
    \Ad_{\theta}(j\ j+k)
    \quad\text{where}\quad
    \forall \ell\left(n+\ell=\ell+1\right).
  \end{equation*}
  Note in particular that the cycle $\eta:=(1\ 2\ \cdots\ n)$ decomposes as $\xi_j\cdots \xi_{j+n-1}$ for any $j\in [n]$, so that
  \begin{equation}\label{eq:te.decomp}
    \forall \left(j\in [n]\right)
    \quad:\quad
    \theta\cdot \eta
    =
    \theta_{jn}\cdots \theta_{j1}\cdot \theta.
  \end{equation}
  Applying \Cref{eq:trnsp.lipsch} repeatedly to each \Cref{eq:te.decomp} shows that the single value $\tau_{\theta\eta}$ is expressible, for each $j\in [n]$, as 
  \begin{equation*}
    \tau_{\theta\eta}
    =
    \theta_{\mathbf{i}_j}\cdot \tau_{\theta}
    :=
    \prod^{\text{ordered}}_{i\in \mathbf{i}_j}\theta_{ji}
    \cdot \tau_{\theta}
  \end{equation*}
  for (possibly empty) tuples
  \begin{equation*}
    \mathbf{i}_j
    :=
    \left(i_{jk}>\cdots>i_{j1}\right)
    \subseteq [n]
  \end{equation*}
  (note that the $k$s also depend on $j$; that dependence is suppressed for legibility). Observe next that we $\theta_{\mathbf{i}_j}$ cannot all be equal for $j\in [n]$ unless
  \begin{enumerate}[(I),wide]
  \item\label{item:pr:chartaus:fl} all $\mathbf{i}_j$ are full;
  \item\label{item:pr:chartaus:mpty} or all are empty. 
  \end{enumerate}
  Indeed, a non-empty product $(j\ j+k_{\alpha})\cdots (j\ j+k_1)$ with the set $\{k_i\}$ avoiding at least some $s\in \left([n]\setminus \{j\}\right)$ will fix that $s$ while at the same time being non-trivial. An equality
  \begin{equation*}
    (j\ j+k_{\alpha})\cdots (j\ j+k_1)
    =
    (s\ s+\ell_{\beta})\cdots (s\ s+\ell_1)
  \end{equation*}
  would then force the latter product to fix $s$, so that it must be empty. But this contradicts its non-triviality.
  
  Case \Cref{item:pr:chartaus:mpty} means that not only does $\theta$ linger along all c-simple $\xi_j$, but also every $\theta\cdot \xi_j$ lingers along every $\xi_{j+1}$ (and hence along all $\xi$). It follows that everything in sight lingers along every $\xi$, hence \Cref{item:pr:chartaus:lngr} above. Similarly, \Cref{item:pr:chartaus:fl} begets \Cref{item:pr:chartaus:grw}.
\end{proof}


\begin{remark}\label{re:need.cyc}
  It was crucial, in \Cref{pr:chartaus}, that the transpositions $\sigma$ of \Cref{eq:trnsp.lipsch} range over the full contingent of $n$ c-simples: only the generators $(j\ j+1)$, $1\le j\le n-1$ of $S_n$ would not have sufficed, per \Cref{ex:just.n1}. 
\end{remark}

\begin{example}\label{ex:just.n1}
  The following self-map of $S_3$ is neither a right translation nor constant, but nevertheless satisfies \Cref{eq:trnsp.lipsch} with $\sigma$ ranging only over the two generators $(1\ 2)$ and $(2\ 3)$ of $S_3$. 
  \begin{equation*}
    \emptyset,\ (2\ 3)\xmapsto{\quad}\emptyset
    ,\quad
    (1\ 2)\xmapsto{\quad}(1\ 2)
    ,\quad
    (1\ 3),\ (1\ 3\ 2)\xmapsto{\quad}(1\ 3)
    ,\quad
    (1\ 2\ 3)\xmapsto{\quad}(1\ 2\ 3),
  \end{equation*}
  where `$\emptyset$' stands for the identity of $S_n$ in order to avoid confusion between it and the symbol `1' in cycle-decomposition notation.
\end{example}

\Cref{le:dichot} below is a preliminary remark moving us closer to eventually proving \Cref{thn:cls.sun}. Some terminology will help streamline the statement.

\begin{definition}\label{def:dichot.diag}
  Let $\T\le \SU(n)$ be a maximal torus and $\T\xrightarrow{\phi}M_n$ a continuous map preserving both spectra and commutativity, falling into one of two qualitatively distinct types by \Cref{th:su.tor}. 

  We say that $\phi$ \emph{reorders} (or \emph{is a reordering}) if it is of type \Cref{item:th:su.tor:ord} and \emph{conjugates} (or simply \emph{is a conjugation}) if it is of type \Cref{item:th:su.tor:cj} instead. 
\end{definition}

\begin{lemma}\label{le:dichot}
  The restrictions of continuous CS-preserving map $\SU(n)\to M_n$ to maximal tori either all conjugate or all reorder in the sense of \Cref{def:dichot.diag}.
\end{lemma}
\begin{proof}
  Maximal tori being mutually conjugate \cite[Theorem IV.1.6]{btd_lie_1995}, they constitute the connected space $\cT\left(\SU(n)\right)\cong \SU(n)/N(\T)$ (with $N(\bullet)$ denoting normalizers) for any fixed maximal torus $\T$. To conclude, observe that
  \begin{equation*}
    \cT(\SU(n))
    =
    \left\{\T\ :\ \phi|_{\T}\text{ conjugates}\right\}
    \sqcup
    \left\{\T\ :\ \phi|_{\T}\text{ reorders}\right\}
  \end{equation*}
  is a disjoint union into closed subsets (so that one must be empty by connectedness). Indeed, the continuity of $\phi$ makes both the conjugation and reordering conditions closed:
  \begin{itemize}[wide]
  \item conjugation can be expressed (for continuous maps already known to preserve commutativity and spectra) as $\phi|_{\T}$ being a group morphism;

  \item while reordering can be phrased as constancy of $\phi|_{\T}$ along Weyl-group orbits.
  \end{itemize}
\end{proof}

Some notation extending that introduced in \Cref{thn:cls.sun} will help with the latter's proof.

\begin{notation}\label{not:eigs}
  Recall the $\lambda_j(X)$ of \Cref{thn:cls.sun}. More generally:
  \begin{equation*}
    \forall\left(S\subseteq [n]\right)
    \forall(X\in \SU(n))
    \quad:\quad
    \lambda_S(X)
    :=
    \text{multiset }
    \left\{\lambda_j(X)\ :\ j\in S\right\}.
  \end{equation*}
  Also:
  \begin{equation*}
    \forall\left(S\subseteq [n]\right)
    \quad:\quad    
    E_S(X):=\sum_{j\in S}\left(\text{$\lambda_j$-eigenspace of $X$}\right).
  \end{equation*}
  Note that $\dim E_S(X)\ge |S|$, with equality whenever $\lambda_S(X)$ and $\lambda_{[n]\setminus S}(X)$ are disjoint, which situation we reference by calling $U$ \emph{$S$-isolated}.

  Similarly, write
  \begin{equation*}
    \forall\left(\Lambda\subseteq \bC\right)
    \forall\left(T\in M_n\right)
    \quad:\quad
    E_{\Lambda}(T):=\sum_{\lambda\in \Lambda}\left(\text{$\lambda$-eigenspace of $T$}\right),
  \end{equation*}
  so that $E_S(X)=E_{\lambda_S(X)}(X)$ for $X\in \SU(n)$. We will frequently omit braces in indicating singleton subscripts, as in $E_j$, $\lambda_j$, etc. 
\end{notation}

We also write
\begin{equation*}
  \bG(d,V)
  :=
  \left\{\text{$d$-dimensional subspaces of }V\right\}
  ,\quad
  \bG(V)\left(\text{or }\bG V\right):=\bigcup_d \bG(d,V)
\end{equation*}
for the various \emph{Grassmannians} (\cite[Example 1.36]{lee2013introduction}, \cite[Example 1.1.3]{ww_twst}) of a finite-dimensional (mostly complex) vector space $V$. 

We record the following observation (cf. its parallel \cite[Lemma 2.3]{2501.06840v2}). 

\begin{lemma}\label{le:def.grs.crsp}
  Let $n\in \bZ_{\ge 1}$ and $\SU(n)\xrightarrow{\phi}M_n$ be a continuous commutativity- and spectrum-preserving map. For an initial or terminal segment
  \begin{equation*}
    S
    =
    \big([n]_{\le k}
      :=
      [k]
      =
      \left\{1..k\right\}
    \big)
    \ \text{or}\ 
    \big(
      [n]_{> k}
      :=
      \left\{k+1..n\right\}
    \big)
  \end{equation*}
  the correspondence
  \begin{equation}\label{eq:psi.corresp}
    \bG(|S|,\bC^n)
    \ni
    W=E_S(U)
    \xmapsto[\quad\text{some $S$-isolated }U\in \SU(n)\quad]{\quad\Psi_{S}=\Psi_{\phi,S}\quad}
    E_{S}(\phi U)
    \in
    \bG(|S|,\bC^n)
  \end{equation}
  is a well-defined continuous map, independent of the choice of $S$-isolated $U\in \SU(n)$ in the sense of \Cref{not:eigs}.
\end{lemma}
\begin{proof}
  The initial/terminal branches being perfectly analogous, we handle the former for $S=[k]$, $1\le k<n$ (there being nothing to prove for $k=n$). Continuity will moreover follow immediately from that of $\phi$, so we focus on $\Psi_S=\Psi_{[k]}$ being well defined: we fix an $S$-isolated $U\in \SU(n)$ with $k$-dimensional $V:=E_S(U)$ and argue that $E_S(\phi U)=E_S(\phi U_V)$ for $S$-isolated $U_V\in \SU(n)$ depending solely on $V$. Thus:
  \begin{itemize}[wide]
  \item select $a<b\le a+1\in \bR$ with $ka+(n-k)b=0$ once for the duration of the proof (so that the choice depends only on $k$);

  \item set
    \begin{equation*}
      U_V
      :=
      \left(\exp(2\pi i a)\text{ on }V\right)
      \oplus
      \left(\exp(2\pi i b)\text{ on }V^{\perp}\right);
    \end{equation*}

  \item and observe that $U$ and $U_V$ commute, so belong to a common maximal torus; the conclusion (that $E_S(\phi U)=E_S(\phi U_V)$) follows from \Cref{le:dichot}.
  \end{itemize}
\end{proof}

For subspaces $V,V'\le \bC^n$ write $V\text{\circled{$\perp$}} V'$ (\emph{weak perpendicularity}) if the orthogonal projections onto $V$ and $V'$ commute. Equivalently:
\begin{equation*}
  V\text{\circled{$\perp$}} V'
  \xLeftrightarrow{\quad}
  \left(V\cap V'\right)
  \perp
  \left(V\ominus \left(V\cap V'\right)\right)
  +
  \left(V'\ominus \left(V\cap V'\right)\right),
\end{equation*}
`$\ominus$' denoting the orthogonal complement of one space in another. 

In preparation for the arguably less interesting (not-quite) half of \Cref{thn:cls.sun}, disposed of easily enough in \Cref{pr:thn:cls.sun:ord}, we need the following remark. 

\begin{lemma}\label{le:if.wperp}
  If a continuous CS preserver $\SU(n)\xrightarrow{\phi}M_n$, $n\in \bZ_{\ge 1}$ reorders on at least one maximal torus $\T\le \SU(n)$ and $\Psi_S=\Psi_{\phi,S}$ of \Cref{le:def.grs.crsp} is well-defined for some $S\subseteq [n]$, then
  \begin{equation*}
    V
    \text{\circled{$\perp$}}
    V'
    \xRightarrow{\quad}
    \Psi_{S}(V)
    =
    \Psi_{S}(V').
  \end{equation*}
\end{lemma}
\begin{proof}  
  Let $T\in \U(n)$ fix $\left(V\cap V'\right)\oplus \left(V+V'\right)^{\perp}$ pointwise and map $V\ominus \left(V\cap V'\right)$ onto $V'\ominus \left(V\cap V'\right)$. For an $S$-isolated $U$ with $E_S(U)=V$ we have
  \begin{equation*}
    E_S\left(U':=\Ad_T U\right)
    =
    V'
    \quad\text{and}\quad
    U,U'\text{ commute},
  \end{equation*}
  so that
  \begin{equation*}
    \Psi_S(V)
    =
    E_S(\phi U)
    \xlongequal{\ \text{\Cref{le:dichot}}\ }
    E_S(\phi U')
    =
    \Psi_S(V')
  \end{equation*}
  by applying reordering on any maximal torus containing $U,U'$. 
\end{proof}

\begin{corollary}\label{cor:psi.ct}
  If a continuous CS preserver $\SU(n)\xrightarrow{\phi}M_n$, $n\in \bZ_{\ge 3}$ reorders on at least one maximal torus $\T\le \SU(n)$ any well-defined $\Psi_S=\Psi_{\phi,S}$ in \Cref{le:def.grs.crsp} is constant. 
\end{corollary}
\begin{proof}
  Let $1\le k:=|S|<n$. The conclusion will follow from \Cref{le:if.wperp} once we argue that because (crucially: \Cref{ex:n2.ord.ncnst}) $n\ge 3$, any two $k$-dimensional subspaces $V,V'\le \bC^n$ can be linked through a finite chain
  \begin{equation*}
    V=:V_0
    \text{\circled{$\perp$}}
    V_1
    \text{\circled{$\perp$}}
    \cdots
    \text{\circled{$\perp$}}
    V_s
    :=
    V'.
  \end{equation*}
  If $k\le n-2$, any one line $\ell\le V=\ell\oplus W$ can be exchanged for another, $\ell'$, via
  \begin{equation*}
    \left(\ell\oplus W\right)
    \text{\circled{$\perp$}}
    \left(\ell''\oplus W\right)
    \text{\circled{$\perp$}}
    \left(\ell'\oplus W\right)
    \quad\text{for}\quad
    \ell''\perp \left(\ell+\ell'+W\right),
  \end{equation*}
  permitting the gradual substitution of $V'$ for $V$.
  
  On the other hand, for $k=n-1$ (and $V\ne V'$) set
  \begin{equation*}
    \ell:=V\ominus\left(W:=V\cap V'\right)
    ,\quad
    \ell':=V'\ominus W,
  \end{equation*}
  pick an arbitrary line $\ell''\in W$, and observe that $\left(W\ominus \ell''\right)\oplus\ell\oplus \ell'$ is weakly-orthogonal to both $V$ and $V'$. 
\end{proof}

\begin{proposition}\label{pr:thn:cls.sun:ord}
  Let $n\in \bZ_{\ge 3}$ and $\SU(n)\xrightarrow{\phi}M_n$ be a continuous CS-preserving map.

  If $\phi|_{\T}$ reorders for at least one maximal torus $\T\le \SU(n)$ then it is of the type listed as \Cref{item:thn:cls.sun:ord} in \Cref{thn:cls.sun}. 
\end{proposition}
\begin{proof}

  We again have reordering on \emph{all} maximal tori, per \Cref{le:dichot}. We know from \Cref{le:def.grs.crsp} and \Cref{cor:psi.ct} that all initial- or terminal-segment $S\subseteq [n]$ yield constant $\Psi_S$. Given that
  \begin{equation*}
    \forall\left(j\in [n]\right)
    \left(
      \im \Psi_{j}
      =
      \im\Psi_{[n]\le j}
      \cap
      \im\Psi_{[n]\ge j}
    \right),
  \end{equation*}
  so too are all singleton-indexed $\Psi_j$, $j\in [n]$. This suffices to conclude. 
\end{proof}

\Cref{pr:thn:cls.sun:ord} certainly does \emph{not} hold for $n=2$:

\begin{example}\label{ex:n2.ord.ncnst}
  For any 
  \begin{equation*}
    \bP \bC^2:=\bG(1,\bC^2)
    \ni \ell
    \xrightarrow[\quad\text{continuous}\quad]{\quad\omega\quad}
    \PGL(2):=\GL(2)/\text{scalars}
    ,\quad
    \omega(\ell)=\omega(\ell^{\perp})
  \end{equation*}
  the map
  \begin{equation*}
    \SU(2)
    \ni
    X
    \xmapsto{\quad}
    \Ad_{\omega(E_1(X))}\left(\lambda_1(X),\ \lambda_2(X)\right)
    \in M_2
  \end{equation*}
  is continuous and CS-preserving, and reorders on each maximal torus but not ``globally'' (so is neither of the form \Cref{item:thn:cls.sun:cj} nor \Cref{item:thn:cls.sun:ord} in the language of \Cref{thn:cls.sun}) provided $\omega$ is non-constant modulo the invertible diagonal matrices. 
\end{example}

\pf{thn:cls.sun}
\begin{thn:cls.sun}
  That maps of type either \Cref{item:thn:cls.sun:cj} or \Cref{item:thn:cls.sun:ord} are continuous and CS-preserving is self-evident, so it is the converse that we are concerned with. At this stage we know that a continuous CS-preserving map $\SU(n)\xrightarrow{\phi}M_n$
  \begin{itemize}[wide]
  \item restricts to every maximal torus as either a conjugation or a reordering (by \Cref{th:su.tor});

  \item so must be of type \Cref{item:thn:cls.sun:ord} if reordering on at least one maximal torus, by \Cref{pr:thn:cls.sun:ord}. 
  \end{itemize}
  What it remains to argue, then, is that if $\phi$ conjugates on \emph{every} maximal torus then it must be of type \Cref{item:thn:cls.sun:cj}. We can now simply outsource the conclusion to the unitary (as opposed to \emph{special} unitary) analogue \cite[Theorem 2.1]{2501.06840v2} of \Cref{thn:cls.sun}.

  Observe first that the conjugation-on-tori assumption implies the scaling compatibility of $\phi$:
  \begin{equation}\label{eq:zx}
    \forall\left(\zeta\in \bS^1\cap \SU(n)\cong \bZ/n\right)
    \forall\left(X\in \SU(n)\right)
    \quad:\quad
    \phi(\zeta X)=\zeta \phi(X). 
  \end{equation}
  This suffices to ensure that $\phi$ admits a continuous, CS-preserving extension to all of $\U(n)$: take \Cref{eq:zx} as the \emph{definition} of that extension, allowing $\zeta$ to range over the entire central circle $\bS^1\le \U(n)$. The conclusion now follows from the aforementioned \cite[Theorem 2.1]{2501.06840v2}, which says (among other things) that continuous CS-preserving maps $\U(n)\to M_n$ are of type \Cref{item:thn:cls.sun:cj}. 
\end{thn:cls.sun}

\begin{remarks}\label{res:post.sun.cls}
  \begin{enumerate}[(1),wide]
  \item\label{item:res:post.sun.cls:ord.no.scl} It was essential, in the proof just given, that we dispose of \Cref{item:thn:cls.sun:ord} before extending $\phi$ to all of $\U(n)$ by scaling: reordering maps (i.e. those of type \Cref{item:thn:cls.sun:ord}) are constant along conjugacy classes (for they depend only on the spectra of their arguments), so cannot satisfy \Cref{eq:zx}.

    If, say, for some $n^{th}$ root of unity $\zeta$ the operators $\zeta X$ and $X$ are mutual conjugates (e.g. $X=\left(\zeta^j\right)_{j=0}^{n-1}$ for primitive $\zeta^n=1$ and odd $n$), then $\phi(\zeta X)=\phi(X)$ for the maps $\phi$ of \Cref{thn:cls.sun}\Cref{item:thn:cls.sun:ord}.

  \item\label{item:res:post.sun.cls:rev.strtg} It is perhaps apposite at this point to note that the proof strategy for \Cref{thn:cls.sun} can be reversed: its branch \Cref{item:thn:cls.sun:cj} can be treated very much along the lines of the unitary version of \cite[Theorem 2.1]{2501.06840v2}:
    \begin{itemize}[wide]
    \item One would start the proof as before, by setting aside the reordering case \Cref{item:thn:cls.sun:ord} and assuming throughout the proof that the continuous CS-preserving map $\phi$ conjugates along all maximal tori.
      
    \item In that case, the continuous self-map $\Psi_1=\Psi_{\phi,1}$ of $\bP^1:=\bG(1,\bC^n)$ introduced in \Cref{le:def.grs.crsp} meets the hypotheses of the \emph{Fundamental Theorem of Projective Geometry} \cite[Theorem 3.1]{zbMATH01747827} so (as in \cite[Proposition 2.5]{2501.06840v2}) we have
      \begin{equation*}
        \bP^1\ni \ell
        \xmapsto{\quad\Psi_1\quad}
        T\ell
        \in \bP^1
      \end{equation*}
      for a linear or conjugate-linear invertible $T$ on $\bC^n$.

    \item Then, as in the proof of \cite[proof of Theorem 2.1, unitary case]{2501.06840v2}, this gives the desired description for $\phi$: conjugation by $T$ if the latter is linear, and $\Ad_{TJ}(-)^t$ if $T$ is conjugate-linear for an appropriately-chosen (also conjugate-linear) $J$. 
    \end{itemize}
    That proof in hand, one could then recover the unitary version of \cite[Theorem 2.1]{2501.06840v2} from \Cref{thn:cls.sun} (rather than the other way round): see \Cref{pr:su2u} below.
  \end{enumerate}  
\end{remarks}

\begin{proposition}\label{pr:su2u}
  Assuming \Cref{thn:cls.sun}, every continuous CS-preserving map $\U(n)\to M_n$, $n\in \bZ_{\ge 3}$ is of type \Cref{item:thn:cls.sun:cj}.
\end{proposition}
\begin{proof}
  Observe first that continuous CS-preserving maps $\U(n)\xrightarrow{\phi} M_n$ must be homogeneous (i.e. intertwine scalars): one can either invoke \cite[Proposition 2.2]{2501.06840v2} or simply note that for every maximal torus $\T\le \U(n)$ $\phi$ restricts to a conjugation on every connected component of
  \begin{equation*}
    \left\{\text{simple-spectrum unitaries in }\T\right\}\subseteq \T
  \end{equation*}
  and such connected components are invariant under the connected scalar subgroup $\bS^1\le \U(n)$. 
  
  Because $\phi$ restricts to a continuous CS-preserving map on $\SU(n)$, the conclusion follows from \Cref{thn:cls.sun} after noting that the reordering-type maps of \Cref{thn:cls.sun}\Cref{item:thn:cls.sun:ord} are not homogeneous (as pointed out in \Cref{res:post.sun.cls}\Cref{item:res:post.sun.cls:ord.no.scl}) and hence do not extend to $\U(n)$.
\end{proof}

\pf{thn:cls.sa}
\begin{thn:cls.sa}
  The argument in the first part of the proof of \Cref{pr:su2u}, delivering the homogeneity of a continuous CS-preserving map on $\U(n)$, functions also to show that any such map $H_n\xrightarrow{\phi}M_n$ intertwines affine transformations:
  \begin{equation*}
    \phi(\alpha X+\beta) = \alpha\phi(X)+\beta
    ,\quad
    \forall X\in H_n,\ \alpha\in \bC^{\times}
    \ \text{and}\ 
    \beta\in \bC
  \end{equation*}
  (cf. \cite[Corollary 4]{Petek-HM}). In particular, if and when convenient, it suffices to prove the conclusion for the restriction of $\phi$ to
  \begin{equation*}
    H^{\le}_n
    :=
    \left\{
      X\in H_n
      \ :\
      \mathrm{trace}~X=0
      \ \wedge\
      \eta_n(X)-\eta_1(X)\le 1
    \right\}
  \end{equation*}
  with $\eta_i$ as in \Cref{eq:sa.reord.spc}. Per the discussion following \Cref{thn:cls.sun}, \cite[Lemmas 1 and 2]{MR210096} imply that
  \begin{equation*}
    H_n
    \xrightarrowdbl{\quad\exp(2\pi i\cdot)}
    \U(n)
  \end{equation*}
  \emph{almost} restricts to a homeomorphism $H^{\le}_n\xrightarrow{\Theta}\SU(n)$: onto, and identifying $X_0\ne X_1$ precisely in the boundary cases when
  \begin{align*}
    \left\{\eta_1(X_0),\ \eta_n(X_0)\right\}
    =
    \left\{\eta_1(X_1),\ \eta_n(X_1)\right\}
    ,\quad
    \eta_n\left(X_{0,1}\right)-\eta_1\left(X_{0,1}\right)=1
  \end{align*}
  and
  \begin{equation*}
    \left(\text{$\eta_1$-eigenspace of $X_{i}$}\right)
    \oplus
    \left(\text{$\eta_n$-eigenspace of $X_i$}\right)
  \end{equation*}
  is independent of $i=0,1$. \Cref{thn:cls.sun} will apply to 
  \begin{equation*}
    \begin{tikzpicture}[>=stealth,auto,baseline=(current  bounding  box.center)]
      \path[anchor=base] 
      (0,0) node (l) {$\SU(n)$}
      +(2,.5) node (ul) {$H^{\le}_n$}
      +(4,.5) node (ur) {$M_n$}
      +(6,0) node (r) {$M_n$}
      ;
      \draw[->] (l) to[bend left=6] node[pos=.5,auto] {$\scriptstyle \Theta^{-1}$} (ul);
      \draw[->] (ul) to[bend left=6] node[pos=.5,auto] {$\scriptstyle \phi$} (ur);
      \draw[->] (ur) to[bend left=6] node[pos=.5,auto] {$\scriptstyle \exp(2\pi i\cdot)$} (r);
      \draw[->] (l) to[bend right=6] node[pos=.5,auto] {$\scriptstyle $} (r);
    \end{tikzpicture}
  \end{equation*}
  to deliver the conclusion (for $H^{\le}_n$, hence also $H_n$) provided that map is well defined: we are left having to argue that $\phi$ is compatible with $\Theta^{-1}$, in the sense that
  \begin{equation}\label{eq:phi.theta.compat}
    \Theta(X)=\Theta(X')
    \xRightarrow{\quad}
    \Theta(\phi X)=\Theta(\phi X').
  \end{equation}
  Set
  \begin{equation*}
    \begin{aligned}
      X:=\bP\bC^n,\quad
      X_{\perp}^n
      &:=
        \left\{\text{orthogonal tuples of lines in $\bC^n$}\right\}\\
      X_{\spn}^n
      &:=
        \left\{\text{spanning line $n$-tuples}\right\}            
    \end{aligned}
  \end{equation*}
  and for
  \begin{equation*}
    \lambda
    =
    \left(\lambda_1\le \cdots\le \lambda_n\right)
    \in
    \bR^{n,\le}
    :=
    \left\{\text{non-decreasing $n$-tuples}\right\}
    \subseteq \bR^n
  \end{equation*}
   and $\ell=\left(\ell_i\right)_{i=1}^n\in X^n_{\spn}$ denote by $T_{\ell,\lambda}$ the operator scaling each $\ell_i$ by $\lambda_i$. $\phi$ induces a continuous map
  \begin{equation*}
    X_{\perp}^n
    \ni
    \ell
    \xmapsto{\quad}
    \wt{\phi}\ell
    \in
    X_{\spn}^n
    :=
    \left\{\text{spanning line $n$-tuples}\right\}
  \end{equation*}
  defined by
  \begin{equation*}
    \forall \lambda\in \bR^{n,\le}
    \quad:\quad
    H_n
    \ni
    T_{\ell,\lambda}
    \xmapsto{\quad\phi\quad}
    T_{\wt{\phi}\ell,\lambda}:
  \end{equation*}
  the independence on $\lambda$ follows from continuity and commutativity preservation by deforming \emph{simple} (i.e. distinct-entry) tuples
  \begin{equation*}
    \lambda
    =
    \left(\lambda_1<\cdots<\lambda_n\right)
    \in
    \bR^{n,<}
    :=
    \left\{\text{increasing $n$-tuples}\right\}
    \subseteq \bR^n
  \end{equation*}
  into one another continuously and passing to arbitrary $\lambda$ by continuity. 
  
  The permutation $S_n$-action $\triangleright$ on $X^n$ restricts to actions on both $X^n_{\perp}$ and $X^n_{\spn}$, and \Cref{eq:phi.theta.compat} amounts to showing that
  \begin{equation*}
    \phi T_{(1\ n)\triangleright \ell,\lambda}
    \in
    \left\{
      T_{(1\ n)\triangleright \wt{\phi}\ell,\lambda}
      ,\
      T_{\wt{\phi}\ell,\lambda}
    \right\}
    \quad
    \bigg(
    \xLeftrightarrow{\quad}
    \ 
    \wt{\phi}\left((1\ n)\triangleright\ell\right)
    \in
    \left\{
      (1\ n)\triangleright \wt{\phi}\ell
      ,\
      \wt{\phi}\ell
    \right\}
    \bigg).
  \end{equation*}
  For $\lambda\in \bR^{n,<}$ and $j\in \left\{1..n-1\right\}$ one can deform $T_{\ell,\lambda}$ into $T_{(j\ j+1)\triangleright\ell,\lambda}$ through a homotopy interchanging the eigenvalues along $\ell_j$ and $\ell_{j+1}$ and leaving all else unaffected, so that we have 
  \begin{equation*}
    \forall\left(1\le j\le n-1\right)
    \quad:\quad
    \wt{\phi}\left((j\ j+1)\triangleright\ell\right)
    \in
    \left\{
      (j\ j+1)\triangleright \wt{\phi}\ell
      ,\
      \wt{\phi}\ell
    \right\},
  \end{equation*}
  with the choice between the two options independent of the $\ell\in X^n_{\perp}$ by (the continuity of $\phi$ and) the latter space's connectedness. Having just observed that the base case $k=1$ holds, we will argue for the validity of $\cP_{j,j+k}$ for all $1\le j\le n-k$ where
  \begin{equation*}
    \cP_{a,b}
    \quad:\quad
    \wt{\phi}\left((a\ b)\triangleright\ell\right)
    \in
    \left\{
      (a\ b)\triangleright \wt{\phi}\ell
      ,\
      \wt{\phi}\ell
    \right\};
  \end{equation*}
  the argument inducts on $1\le k\le n-1$, applying
  \begin{equation}\label{eq:abbcac}
    \cP_{a,b}\wedge \cP_{b,c}    
    \xRightarrow{\quad}
    \cP_{a,c}
  \end{equation}
  to
  \begin{equation*}
    (a,b,c)
    :=
    \left(j,\ j+k-1,\ j+k\right).
  \end{equation*}
  To confirm \Cref{eq:abbcac}, write $\sigma^{\circ}$ for $\sigma\in S_n$ and 
  \begin{equation*}
    \forall\left(\ell\in X^n_{\perp}\right)
    \left(
      \wt{\phi}\left(\sigma\triangleright \ell\right)
      =
      \sigma^{\circ}\triangleright \wt{\phi}\ell
    \right),
  \end{equation*}
  so that the hypothesis of \Cref{eq:abbcac} reads
  \begin{equation}\label{eq:abbc.circ}
    (a\ b)^{\circ}\in \left\{(a\ b),\ 1\right\}
    \quad\text{and}\quad
    (b\ c)^{\circ}\in \left\{(b\ c),\ 1\right\}.
  \end{equation}
  We then have 
  \begin{equation*}
    \forall\left(\sigma_{1\le i\le m}\in \left\{(a\ b),\ (b\ c)\right\}\right)
    \left(
      \wt{\phi}\left(\sigma_1\cdots \sigma_m\triangleright \ell\right)
      =
      \sigma_1^{\circ}\cdots \sigma_m^{\circ}\triangleright \wt{\phi}\ell
    \right),
  \end{equation*}
  hence
  \begin{equation*}
    (a\ b)(b\ c)(a\ b)
    =
    (a\ c)
    =
    (b\ c)(a\ b)(b\ c)
    \xRightarrow{\quad}
    (a\ b)^{\circ}(b\ c)^{\circ}(a\ b)^{\circ}
    =
    (b\ c)^{\circ}(a\ b)^{\circ}(b\ c)^{\circ}
  \end{equation*}  
  Given \Cref{eq:abbc.circ}, this is only possible if
  \begin{equation*}
    \bigg(
      (a\ b)^{\circ}=(a\ b)\wedge (b\ c)^{\circ}=(b\ c)
    \bigg)
    \vee
    \bigg(
      (a\ b)^{\circ}=1=(b\ c)^{\circ}
    \bigg).
  \end{equation*}
  This confirms $\cP_{a,c}$, completing the proof of \Cref{eq:abbcac} and thus of the theorem. 
\end{thn:cls.sa}

\section{Lipschitz self-maps of Coxeter groups}\label{se:cox.lip}

It might be of some interest to observe that \Cref{pr:chartaus} is an instance of a wider pattern, to be further examined in \Cref{th:cox.lip}: the latter applies to \emph{Coxeter systems} \cite[p.2]{bb_comb-cox} $(W,S)$ and their underlying \emph{Coxeter groups}, $W$, with $S_n$ realized as one such as usual, via \cite[Example 1.2.3]{bb_comb-cox}: the system $S\subseteq S_n$ of generators is
\begin{equation*}
  S
  :=
  \left\{(i\ i+1)\ :\ 1\le i\le n-1\right\}. 
\end{equation*}
For background on Coxeter groups we refer the reader to standard sources such as \cite{bb_comb-cox,hmph_cox}, with more specific citations where needed. The following piece of vocabulary is meant as reminiscent of the \emph{Lipschitz} maps ubiquitous \cite[\S 1.4]{bbi} in metric geometry, providing shorthand for \Cref{eq:trnsp.lipsch}.

\begin{definition}\label{def:s.lip}
  Let $W$ be a group and $W'\xrightarrow{\phi}2^W$ a partial function for $W'\subseteq W$.

  A self-map $\tau_{\bullet}\in W^W$ is \emph{(right-)$\phi$-Lipschitz} if
  \begin{equation}\label{eq:phi.lip}
    \forall\left(\theta\in W'\right)
    \forall\left(\sigma\in \phi(\theta)\right)      
    \quad:\quad
    \tau_{\sigma\theta}\in \left\{\tau_{\theta},\ \sigma\tau_{\theta}\right\}.
  \end{equation}
  When $\phi$ takes a constant value $T\in 2^W$ we refer to $\tau$ as \emph{(right-)$T$-Lipschitz}. Explicitly:
    \begin{equation}\label{eq:tlip}
      \forall\left(\theta\in W\right)
      \forall\left(\sigma\in T\right)      
      \quad:\quad
      \tau_{\sigma\theta}\in \left\{\tau_{\theta},\ \sigma\tau_{\theta}\right\}.
    \end{equation}
    We will mostly be interested in the case when $(W,S)$ is a Coxeter system and $\phi$ takes values in its set $\Ad_W S$ of \emph{reflections} \cite[\S 1.3, p.12]{bb_comb-cox}.
\end{definition}

\begin{remarks}\label{res:runif}
  \begin{enumerate}[(1),wide]
  \item\label{item:res:runif:1lip.brh} It follows from \cite[Theorem 2.2.2]{bb_comb-cox} that for a Coxeter system $(W,S)$ the $S$-Lipschitz condition means precisely that
    \begin{equation*}
      \forall\left(\theta,\eta\in W\right)
      \quad:\quad
      \tau_{\theta\eta^{-1}}\le \theta\eta^{-1}
    \end{equation*}
    for the \emph{Bruhat order} \cite[Deﬁnition 2.1.1]{bb_comb-cox} on $(W,S)$. This is, in other words, the requirement that $\tau$ be \emph{contractive} (i.e. distance non-increasing, or \emph{1-Lipschitz} in the language of \cite[Definition 1.1]{grm_mtr-struct}, say) with respect to a poset-valued distance. 
    
  \item\label{item:res:runif:why.rght} \Cref{def:s.lip} speaks of the \emph{right}-handed Lipschitz property because plainly, condition \Cref{eq:tlip} is invariant under right translation on $W$.
  \end{enumerate}
\end{remarks}

The following observation is immediate.

\begin{lemma}\label{le:mon}
  If $\tau_i\in W^W$, $i=0,1$ are $\phi_i$-Lipschitz respectively, then $\tau_1\circ \tau_0$ is $\phi$-Lipschitz for
  \begin{equation*}
    \left\{w\in \mathrm{dom}~\phi_0\ :\ \phi_0(w)\subseteq \phi_1\left(\tau_w\right)\right\}
    \ni w
    \xmapsto{\quad\phi\quad}
    \phi_0(w).
  \end{equation*}
  In particular, $T$-Lipschitz self-maps constitute a monoid for any $T\subseteq W$.  \qedhere
\end{lemma}



\begin{theorem}\label{th:cox.lip}
  Let $(W,S)$ be a Coxeter system and $T$ its set of reflections. The right $T$-Lipschitz self-maps $\tau\in W^W$ are precisely the constants and right translations. 
\end{theorem}

\begin{remark}\label{re:red.conn}
  The right translations and constants plainly constitute a monoid under composition, contained in that (\Cref{le:mon}) of all $T$-Lipschitz maps. To prove the opposite inclusion it will thus suffice, upon precomposing an arbitrary $T$-Lipschitz map $\tau$ with right translation by $\tau_1^{-1}$, to assume $\tau_1=1$ whenever convenient.

  \Cref{pr:conn.suff} below further restricting the scope to systems $(W,S)$ with connected underlying \emph{Coxeter graphs} \cite[\S 1.1]{bb_comb-cox}, it is in the form of
  \begin{equation}\label{eq:ws.conn}
    \left((W,S)\text{ connected }\ \text{and}\ \tau_1=1\right)
    \quad
    \xRightarrow{\quad}
    \quad
    \tau_{\bullet}\in \left\{\id,\ 1\right\}
  \end{equation}
  that we address the claim after due preparation.  
\end{remark}

\begin{proposition}\label{pr:conn.suff}
  If \Cref{th:cox.lip} holds for connected Coxeter graphs it does in general. 
\end{proposition}
\begin{proof}
  Writing $\left(W_{S'},S'\right)$ for the Coxeter system induced by a set of simple reflections, note that in assuming (\Cref{re:red.conn}) $\tau_1=1$, $\tau$ restricts on every $W_{S'}$ to a self-map satisfying the hypotheses. For $S'\subseteq S$ label
  \begin{itemize}[wide]
  \item by $\cP_{S'\mid 1}$ the statement that (assuming $\tau_1=1$) $\tau_{S'}:=\tau|_{W_{S'}}\equiv 1$;
  \item by $\cP_{S'\mid \id}$ the analogous statement that $\tau_{S'}=\id_{W_{S'}}$;
  \item and $\cP_{S'}:=\cP_{S'\mid 1}\vee \cP_{S',\id}$ (i.e. the claim that the desired conclusion holds for $\tau_{S'}$). 
  \end{itemize}
  $\cP_{S'}$ are assumed valid for connected $S'\subseteq S$ (i.e. sets internally linked in the Coxeter graph), and the goal is to prove all $\cP_{S_0}$ valid by induction on the number of connected components of the induced graph on $S_0$.

  Consider reflections
  \begin{equation*}
    t_{i}\in W_{S_{i}}
    ,\quad
    i=0,1
    ,\quad
    \text{no edges $\emptyset\ne S_0\leftrightarrow S_1\ne \emptyset$}. 
  \end{equation*}
  $\cP_{S_i}$ are assumed valid by induction and for reflections $t_i\in W_{S_i}$ we have
  \begin{equation*}
    \tau_{t_0t_1}\in \left\{1,t_0t_1\right\}
    ,\quad
    \begin{aligned}
      \tau_{t_0t_1}=1&\iff \tau_{t_i}=1,\ i=0,1\\
      \tau_{t_0t_1}=t_{0}t_1&\iff \tau_{t_i}=t_i,\ i=0,1.
    \end{aligned}
  \end{equation*}
  It follows that either $\cP_{S'\mid 1}$ holds for all connected components $S'\subseteq S$, or $\cP_{S'\mid \id}$ does (i.e. there is no admixture: connected components cannot differ in this respect). 

  Assume $\cP_{S'\mid 1}$ for some connected component $S'$ (hence for all of them by the foregoing discussion). To prove that $\tau_w=1$ by induction on the number of connected components involved in a minimal decomposition of $w$, simply observe that for mutually disconnected $S_i\subseteq S$, $i=0,1$ and $w=t_0 t_1$, $1\ne t_i\in S_i$ $\tau_w$ belongs to both $W_{S_i}$.

  If on the other hand $\cP|_{S'\mid \id}$ holds for all connected components, take $t_i\in S_i$, $i=0,1$ as above and decompose $t_i=s_i t'_i$ minimally. Induction on \emph{length} \cite[\S 1.4]{bb_comb-cox} allows the assumption $\tau_{t'_0 t'_1}=t'_0 t'_1$, so that $\tau_{t_0 t_1}\in \left\{t_0 t_1, t'_0 t'_1\right\}$. That the latter possibility cannot obtain follows from $\tau_{t_0}=t_0$, by gradually building $t_0 t_1$ from $t_0$ by multiplication with the simple reflections in $t_1=s_1 t'_1$.
\end{proof}

Short of being trivial, the simplest examples of Coxeter systems are the \emph{dihedral groups} $I_2(m)$, $m\in \bZ_{\ge 2}\sqcup\{\infty\}$ of \cite[Example 1.2.7]{bb_comb-cox}:
\begin{equation*}
  I_2(m):=\braket{r_i,\ i=1,2}/\left(r_i^2=1,\ (r_1r_2)^m=1\right)  
\end{equation*}
with the last relation empty for $m=\infty$. The corresponding Coxeter graph is connected for all $m\ge 3$. 

\begin{proposition}\label{pr:fin.dih}
  \Cref{th:cox.lip} holds for the dihedral groups $I_2(m)$, $m\in \bZ_{\ge 3}\sqcup \{\infty\}$. 
\end{proposition}
\begin{proof}
  We prove the claim in the form of \Cref{eq:ws.conn}, noting that in this case the Coxeter graph \emph{is} connected. 
  
  \begin{enumerate}[(I),wide]
  \item\label{item:pr:fin.dih:pf.loc}\textbf{: For all $w\in W:=I_2(m)$, $\tau_w\in \{1,w\}$.} For reflections $w\in T\subseteq W$ the $T$-Lipschitz property makes this immediate, and we are of course assuming $\tau_1=1$. The non-trivial non-reflections are precisely the \emph{even} elements $s_1\cdots s_{2k}$, $k\ge 1$, $s-I\in S$ of the dihedral group, and every such is a product of two reflections in at least two ways:
    \begin{equation*}
      w=s_1\cdot \left(s_2\cdots s_{2k}\right)
      \quad\text{and}\quad
      w=s'_1\cdot \left(s_1' s_1s_2\cdots s_{2k}\right)
      ,\quad
      s'_1\ne s_1. 
    \end{equation*}
    The Lipschitz property applied to all of those reflections ensures that
    \begin{equation*}
      \tau_w\in
      \left\{1,w,s_1,s_2\cdots s_{2k}\right\}
      \cap
      \left\{1,w,s'_1,s_1's_1s_2\cdots s_{2k}\right\}
      =
      \{1,w\},
    \end{equation*}
    concluding the proof of the claim.
    
  \item\label{item:pr:fin.dih:pf.glob}\textbf{: Conclusion.} On the one hand, if $\tau_w=w$ for a minimal-length expression $w=s_1\cdots s_{\ell}$ then
    \begin{equation*}
      \forall\left(1\le k\le \ell\right)
      \left(\tau_{s_k\cdots s_{\ell}}=s_k\cdots s_{\ell}\right)
    \end{equation*}
    by \Cref{item:pr:fin.dih:pf.loc} the $S$-Lipschitz property and length considerations. On the other hand, if the minimally decomposed $1\ ne w=s_1\cdots s_{\ell}$ is the shortest among $s_k\cdots s_{\ell}$ with $\tau_w=1$, then \Cref{item:pr:fin.dih:pf.loc} and the $S$-Lipschitz property imply $\ell=1$. This suffices to conclude: if $\tau_w=1$ for at least one $w\ne 1$, that constancy obtains globally. 
  \end{enumerate}
\end{proof}

\begin{remark}\label{re:s.lip.only}
  It is not unnatural at this stage to ask whether (vagaries of the proof notwithstanding) \Cref{pr:fin.dih} goes through for $S$- (rather than $T$-)Lipschitz maps. It does not: for any $m\in \bZ_{\ge 3}$ the self-map of $W:=I_2(m)=\Braket{r_1,r_2}$ removing, for every $w\in W$, the terminal $r_2$ letter in a \emph{reduced expression} \cite[\S 1.4]{bb_comb-cox} for $w$ if one such exists will be $S$-Lipschitz, fixing $r_1$ and annihilating $r_2$.

  For $m=3$, say, this is the unique $S$-Lipschitz map acting as
  \begin{equation*}
    w_0=
    r_1r_2r_1=r_2 r_1 r_2
    \xmapsto{\quad}
    r_2 r_1
  \end{equation*}
  on the longest element (its other values are then easily filled in).

  This same gadget functions rather generally: \Cref{ex:fold}.
\end{remark}

\begin{example}\label{ex:fold}
  Consider any Coxeter system $(W,S)$ and declare elements $w,w'\in W$ \emph{$s$-adjacent}, $s\in S$ if $w'=sw$ (reversing the convention of \cite[p.10]{ron_build_2e_2009} for \emph{Coxeter complexes}, in other words). Every reflection $t\in T:=\Ad_WS$ determines a \emph{root} $\alpha_T\ni 1$ \cite[Proposition 2.6]{ron_build_2e_2009}, and the \emph{folding} map $W\to \alpha_T$ of \cite[p.15]{ron_build_2e_2009} is clearly 1-Lipschitz with respect to the Bruhat order and hence (\Cref{res:runif}\Cref{item:res:runif:1lip.brh}) $S$-Lipschitz.

  A more concrete description of such a map would be
  \begin{equation*}
    \tau_w
    :=
    \begin{cases}
      ws&\text{if $\exists$ reduced $w=s_1\cdots s_{k-1} s$}\\
      w&\text{otherwise}
    \end{cases}
    ;
  \end{equation*}
  that this is precisely the folding corresponding to $t:=s\in S$ follows from \cite[Corollary 1.4.6]{bb_comb-cox} and the characterization of roots given in \cite[Proposition 2.6(ii)]{ron_build_2e_2009}.
\end{example}

\pf{th:cox.lip}
\begin{th:cox.lip}
  Per \Cref{re:red.conn}, assume that $(W,S)$ connected and $\tau_1=1$. 

  \begin{enumerate}[(I),wide]

  \item\label{item:th:cox.lip:tau.on.s} \textbf{: $\tau|_S\in \left\{\id_S,\ 1\right\}.$} Or: either $\tau_s=s$ for all $s\in S$, or $\tau_s=1$ for all $s\in S$. This follows from \Cref{pr:fin.dih}: any two generators $s,s'\in S$, if connected in the Coxeter graph of $(W,S)$, generate a finite dihedral group $I_2(m)$, $m\in \bZ_{\ge 3}$. This means that they are simultaneously left invariant or sent to $1$ by $\tau_{\bullet}$. Because we are assuming Coxeter-graph connectedness, any two generators can be linked by a path.    

  \item\label{item:th:cox.lip:ww.id.s} \textbf{: If $\tau_w=w$ for some $w\ne 1$ then $\tau|_S=\id|_S$.} If $w=s_1\cdots s_k$ is a reduced expression for $w$, then $k\ge 1$ because $w\ne 1$. The $S$-Lipschitz property implies that
    \begin{equation*}
      \tau_{s_i\cdots s_k}=s_i\cdots s_k
      ,\quad
      \forall 1\le i\le k,
    \end{equation*}
    so in particular $\tau_{s_k}=s_k$. We then have $\tau_s=s$ for all $s\in S$ by step \Cref{item:th:cox.lip:tau.on.s}. 

  \item\label{item:th:cox.lip:if.id.then.id} \textbf{: If $\tau_y=y$ for some $y\ne 1$ then $\tau=\id$.} We already know from \Cref{item:th:cox.lip:ww.id.s} that $\tau|_S=\id_S$. Let $k\ge 2$ be the minimal length of an element with $\tau_w\ne w$ and $w=s_1\cdots s_k$ a reduced expression for such a word. We then have
    \begin{equation*}
      \tau_w=x:=s_2\cdots s_k,
    \end{equation*}
    and switch focus to the (again $T$-Lipschitz) map
    \begin{equation*}
      \tau'_{\bullet}:= (\cdot x^{-1})\circ\tau_{\bullet}\circ (\cdot x).
    \end{equation*}
    We have $\tau_x=x$ and hence $\tau'_1=1$, and also $\tau'_{s_1}=1$. Step \Cref{item:th:cox.lip:ww.id.s} applied to $\tau'$ yields $\tau'_{s}=1$ for \emph{all} $s\in S$, and in particular for $s_2$. But then
    \begin{equation*}
      \tau'_{s_2}=1
      \xRightarrow{\quad}
      \tau_{s_3\cdots s_k}=\tau_{s_2s_2s_3\cdots s_k}=\tau'_{s_2}\cdot s_2\cdots s_k=s_2\cdots s_k,
    \end{equation*}
    contradicting the minimality of $k$ for the length of an element on which $\tau_{\bullet}$ is not identical. The contradiction proves that there are no $w$ with $\tau_w\ne w$, and we are done.

  \item \textbf{: If $\tau_y=1$ for some $y\ne 1$ then $\tau=1$.} We at least have $\tau|_S=1$ by the preceding step \Cref{item:th:cox.lip:if.id.then.id}, so we again proceed by induction, in similar fashion: suppose $w=s_1\cdots s_k$ is a minimal-length element on which $\tau$ is not $1$, so that by the $S$-Lipschitz property we have $\tau_w=s_1$. This time set
    \begin{equation*}
      \tau'_{\bullet}:= \tau_{\bullet}\circ (\cdot s_2\cdots s_k)
      \quad
      \left(\text{again $T$-Lipschitz}\right).
    \end{equation*}

    The sequel is much as before: $\tau'_1=1$ and $\tau'_{s_1}=s_1$, so that $\tau'=\id$ by \Cref{item:th:cox.lip:if.id.then.id}; this contradicts
    \begin{equation*}
      \tau'_{s_2}=\tau_{s_2s_2s_3\cdots s_k}=\tau_{s_3\cdots s_k}=1
    \end{equation*}
    and concludes the proof.  \qedhere
  \end{enumerate}
\end{th:cox.lip}








\addcontentsline{toc}{section}{References}

\def\polhk#1{\setbox0=\hbox{#1}{\ooalign{\hidewidth
  \lower1.5ex\hbox{`}\hidewidth\crcr\unhbox0}}}


\Addresses

\end{document}